\documentclass[12pt]{amsart}

\newtheorem{theorem}{Theorem}[section]
\newtheorem{lemma}[theorem]{Lemma}
\newtheorem{corollary}[theorem]{Corollary}
\theoremstyle{definition}   
\newtheorem{definition}{Definition}

\theoremstyle{remark}

\numberwithin{equation}{section}

\usepackage{times}
\usepackage{enumerate}
\usepackage{tfrupee}
\usepackage{tikz}
\usetikzlibrary{chains,fit}
\usetikzlibrary{shapes,snakes}
\usetikzlibrary{graphs}
\usepackage{graphicx,adjustbox}
\usepackage{hyperref}
\usepackage{amsmath,amssymb}

\newcommand{\LT}{{\rm \mbox{Lt}}}

\usepackage{amscd}
\usepackage{graphicx}
\usepackage[all]{xy}

\title[$d$-sequence and Regular sequence of Quadrics]{$d$-sequence and Regular sequence of Quadrics}
\author{
Joydip Saha
\and
Indranath Sengupta
\and
Gaurab Tripathi
}
\date{}

\address{\small \rm  Stat-Math Unit, Indian Statistical Institute, 203 B.T. Road, Kolkata 700 108.}
\email{saha.joydip56@gmail.com}
\thanks{The author is supported by the NPDF fellowship PDF/2019/001074, sponsored by the SERB, Government of India.}

\address{\small \rm  Discipline of Mathematics, IIT Gandhinagar, Palaj, Gandhinagar, 
Gujarat 382355, INDIA.}
\email{indranathsg@iitgn.ac.in}
\thanks{The second author is the corresponding author; supported by the 
MATRICS research grant MTR/2018/000420, sponsored by the SERB, Government of India.}

\address{\small\rm Department of Mathematics, St.Xavier's College, 
30 Mother Teresa Sarani, Kolkata 700016.}
\email{gelatinx@gmail.com}

\date{}



\begin{document}
\begin{abstract}
Let $K$ be a field and $X$, $Y$ denote matrices such that, the entries 
of $X$ are either indeterminates over $K$ or $0$ and the entries of $Y$ 
are indeterminates over $K$ which are different from those appearing in 
$X$. We consider ideals of the form $I_{1}(XY)$, which is the ideal 
generated by the homogeneous polynomials of degree $2$ given by the 
$1\times 1$ minors of the matrix $XY$. We prove that 
$d$-sequences and regular sequences arise naturally as part of generators 
of $I_{1}(XY)$ for some special cases. We use this information to calculate 
the equations defining the Rees algebra of $I_{1}(XY)$.  
\end{abstract}

\maketitle

\section{Introduction}
Our study originated from the 1974 paper of J. Herzog \cite{herzog} on 
the following theme: Let $R$ be a Noetherian commutative ring with identity. Let 
$\underline{x} = \{x_{1}, \ldots , x_{n}\}$ be a sequence in $R$. Let 
$\mathfrak{a} = (\alpha_{ij})$ be an $m\times n$ matrix with entries in $R$; 
with $m\leq n$. A complex $D_{*}(\underline{x}, \mathfrak{a})$ was constructed 
in \cite{herzog}. Acyclicity conditions on the complex $D_{*}(\underline{x}, \mathfrak{a})$ 
were derived in order to decide perfectness and the Gorenstein property for the 
ideals $\langle a_{1}, \ldots , a_{n}, \Delta \rangle$ and 
$\langle a_{1}, \ldots , a_{n}, \Delta_{1}, \ldots, \Delta_{n} \rangle$, 
where $a_{i} = \sum_{j=1}^{n}\alpha_{ij}x_{j}$, $i=1, \ldots , n$, 
$\Delta = \det(\alpha_{ij})$ (when $n\geq 2$ and $m=n$) and $\Delta_{j}$ is the 
determinant of the matrix obtained from $\mathfrak{a}$ by deleting the $j$-th column (when 
$n\geq 3$ and $m=n-1$). Our aim in this paper is to study a class of ideals of the form 
$I_{1}(XY)$ defined through determinantal conditions, which is similar to the aforesaid class 
studied by Herzog. We show that $d$-sequences 
and regular sequences occur very naturally in the setting of the ideals $I_{1}(XY)$, 
which is the principal object of study in this paper. We show how $d$-sequences occur 
naturally in Theorems \ref{gend} and \ref{altdseq}. In Theorem \ref{maintheorem} we show 
the occurrence of a regular sequence which follows an interesting pattern. However, 
this regular sequence is not a maximal one and the question finding one of maximal 
length remains open. 
\medskip

Ideals of the form $I_{1}(XY)$ have appeared in another significant work 
by Huneke and Ulrich \cite{hu}. We will use some of their observations for proving some 
results on $d$-sequences when $X$ is alternating, in Theorem \ref{altdseq}. In this process 
we would also prove the conjectures proposed in the paper \cite{sstquad}.
\medskip

A sequence of (homogeneous) polynomials $p_{1}, p_{2}, \ldots , p_{r}$ in a 
polynomial ring $R$ is called a \textit{regular sequence} if the ideal 
$\langle p_{1}, p_{2}, \ldots , p_{r} \rangle\neq R$ and if each $p_{i}$ 
is non-zero divisor in $R/\langle p_{1}, p_{2}, \ldots , p_{i-1} \rangle$, 
for every $1\leq i\leq r$. The notion of a $d$-sequence was defined by 
Huneke (see def \ref{dseq}), which generalizes the notion of a regular sequence. 
\medskip

Let Given an ideal $I = \langle a_{1}, \ldots , a_{n}\rangle$ in a Noetherian ring $R$, 
the Rees algebra of $I$ is the $R$-algebra $R[It]= R[a_{1}t, \ldots , a_{n}t]$. One can 
define an $R$-algebra  homomorphism $\pi: R[t_{1}, \ldots , t_{n}] \rightarrow R[It]$ as 
$\pi(t_{i}) = a_{i}t$. The map $\pi$ is a graded map of degree $0$. Therefore $\ker(\pi)$ 
is generated by homogeneous polynomials in $t_{1}, \ldots , t_{n}$. Generators of $\ker(\pi)$ 
are called the equations defining the Rees algebra $R[It]$. Note that the linear polynomials 
$f_{ij} = a_{i}t_{j} - a_{j}t_{i}$ belong to $\ker(\pi)$ and we say that the ideal $I$ is 
of linear type (a notion introduced by Valla) if $\ker(\varphi)$ is generated by the linear 
homogeneous polynomials $f_{ij}$. It was proved independently by C. Huneke and G. Valla, 
and later generalized by K.N. Raghavan, that, if an ideal is generated by a 
$d$-sequence then the ideal is of \textit{linear type}. 
\bigskip

\section{The ideal $I_{1}(XY)$}
Let $R=K[x_{ij},y_{ij} \mid 1\leq i,j\leq n]$, $n\geq 2$, and $x_{ij}$ and $y_{ij}$ are 
distinct indeterminates over the field $K$. Let $X=(x_{ij})_{n\times n}$ and 
$Y=(y_{ij})_{n\times n}$ be generic matrices. Let us write the product of the matrices $X$ and $Y$ as 
$XY=(f_{ij})_{n\times n}$, so that $f_{ij}=\sum_{k=1}^{n}x_{ik}y_{kj}$. 
Let $I_{1}(XY)$ denote the ideal generated by the polynomials $f_{ij}$, 
which are the $1\times 1$ minors of the matrix $XY$. Certain properties 
like primality, primary decomposition and minimal free resolutions have been studied in 
\cite{sstgrob}, \cite{sstprime}, \cite{sstsum}. 
\medskip

It is easy to see that all the 
$f_{ij}$'s defined above do not form a regular sequence. 
For example, if $n=2$ then $x_{12}y_{21}f_{11}+x_{11}y_{12}f_{12}-x_{22}y_{21}f_{21}=x_{21}y_{12}f_{22}$ shows that $f_{11},f_{12},f_{21},f_{22}$ is not a regular sequence. 
However, certain interesting and useful results proved in this direction are the following:
\medskip

\begin{theorem}\label{gensymn}
Let $S = K[x_{ij}, \, y_{j} \mid 1\leq i,j\leq n]$ 
denote the polynomial $K$-algebra and $T_{1},\ldots,T_{n}$ are indeterminates over $S$. Suppose that 
$X =(x_{ij})_{n\times n}$ is either generic or generic symmetric and $Y=(y_{j})_{n\times 1}$ is generic. 
Let $\mathcal{I} = I_{1}(XY) = \langle g_{1},\ldots ,g_{n}\rangle$, where 
$g_{i}=\sum_{j=1}^{n} x_{ij}y_{j}$. The set $\lbrace g_{1},\cdots ,g_{n}\rbrace$ forms a regular sequence. 
Hence, the defining equations of the Rees algebra of $\mathcal{I}$ are only the Koszul relations 
$T_{i}g_{j}-T_{j}g_{i}$, for $1\leq i<j\leq n$.
\end{theorem}

\proof See Theorem 6.1 in \cite{sstgrob}.\qed
\medskip

\begin{theorem}\label{prime}
Let $S = K[x_{ij}, \, y_{j} \mid 1\leq i,j\leq n]$ 
denote the polynomial $K$-algebra and $X=(x_{ij})_{n\times n}$ 
and $Y=(y_{j})_{n\times 1}$ be generic matrices. Suppose that 
$\mathcal{I} = I_{1}(XY) = \langle g_{1},\ldots ,g_{n}\rangle$, where 
$g_{i}=\sum_{j=1}^{n} x_{ij}y_{j}$. Let us define the ideals 
$I_{i}=\langle g_{1},\ldots,g_{i}\rangle$, for $1\leq i\leq n-1$, 
and $I_{n}=\langle g_{1},\ldots,g_{n},\Delta\rangle$, where $\Delta=\det(X)$. 
The ideals $I_{1}, \ldots, I_{n-1}$ and $I_{n}$ are all prime 
ideals.
\end{theorem}

\proof See Corollary 3.2, Theorem 3.3 and Theorem 4.4 in\cite{sstprime}. \qed
\bigskip

\section{$d$-sequence}
\begin{definition}\label{dseq}[Definition 5.5.2; \cite{hunekeswanson}]
Let $A$ be a commutative ring. Set $a_{0}=0$. A sequence of elements $a_{1},\ldots,a_{n}$ 
is said to be a $d$-sequence if 
$$(\langle a_{0},\ldots,a_{i}\rangle:a_{i+1}a_{j})=(\langle a_{0},\ldots,a_{i}\rangle :a_{j}) \quad 
\forall \, 0\leq i\leq n-1, \quad \forall \, j\geq i+1.$$ 
\end{definition}
\medskip

Let $S = K[x_{ij}, \, y_{j} \mid 1\leq i\leq n+1,1\leq j\leq n]$ denote the polynomial 
$K$-algebra and $\widehat{X} = (x_{ij})_{(n+1)\times n}$, $Y=(y_{j})_{n\times 1}$ be 
generic matrices. Suppose $\mathcal{J} = I_{1}(\widehat{X}Y) = \langle g_{1},\ldots ,g_{n+1}\rangle$, where $g_{i}=\sum_{j=1}^{n} x_{ij}y_{j}$. 
We now prove Theorem \ref{gend} to show that $g_{1},\ldots ,g_{n+1}$ is a $d$-sequence. 
This would make the ideal an ideal of 
linear type and the equations defining its Rees algebra are all linear. In other 
words, the Rees algebra and the Symmetric algebra of the ideal are isomorphic. 
We prove the following lemma first:
\medskip

\begin{lemma}\label{cofac} Let $S$, $\widehat{X}$, $Y$ and $I_{1}(\widehat{X}Y)$ be as above. 
Let $\Delta_{i}$ denote the determinant 
of the matrix obtained by deleting $i$'th row of $\widehat{X}$ for $1\leq i\leq n+1$. Then,

\begin{enumerate}[(i)]
\item $\displaystyle\sum_{i=1}^{n+1}(-1)^{i}x_{ij}\Delta_{i}=0$,\, for all $j\in\{1,\ldots ,n\}$.
\item $\displaystyle\sum_{i=1}^{n+1}(-1)^{i}g_{i}\Delta_{i}=0$.
\item $(\langle g_{1},\ldots,g_{n}\rangle :g_{n+1})
=\langle g_{1},\ldots,g_{n},\Delta_{n+1}\rangle$.
\item $(\langle g_{1},\ldots,g_{n}\rangle :\Delta_{n+1})=\langle y_{1},\ldots,y_{n}\rangle$.
\end{enumerate} 
\end{lemma}

\proof (i) \, For each $1\leq j\leq n$, we consider the $(n+1)\times (n+1)$ matrix,  
$$\widetilde{X}_{j}\quad=\quad \begin{pmatrix}
 x_{1j}& x_{11} & x_{12} & \cdots & x_{1n}\\
x_{2j} & x_{21} & x_{22} & \cdots & x_{2n}\\
\vdots & \vdots & \vdots & \vdots & \vdots\\
x_{(n+1)j} & x_{(n+1)1} & x_{(n+1)2} & \cdots & x_{(n+1)(n+1)}\\ 
\end{pmatrix}.$$ 
Then we have $\det(\widetilde{X}_{j})=0$ for all $1\leq j\leq n$. Expanding 
with respect to the first column, we get $\displaystyle\sum_{i=1}^{n+1}(-1)^{i}x_{ij}\Delta_{i}=0$,\, for all $j\in\{1,\ldots ,n\}$. 
\medskip

\noindent (ii)\, Rearranging terms we get, 
$\displaystyle\sum_{i=1}^{n+1}(-1)^{i}g_{i}\Delta_{i}= \displaystyle\sum_{i=1}^{n}(\displaystyle\sum_{i=1}^{n+1}(-1)^{i}x_{ij}\Delta_{i})y_{j}$. 
Therefore $\displaystyle\sum_{i=1}^{n+1}(-1)^{i}g_{i}\Delta_{i}=0$.
\medskip

\noindent (iii) follows from Lemma 6.7 in \cite{sstgrob}.
\medskip

\noindent (iv) By lemma 4.3 in \cite{sstprime}, we have  $\langle y_{1},\ldots,y_{n}\rangle\subset  (\langle g_{1},\ldots,g_{n}\rangle :\Delta_{n+1})$. Let $h\Delta_{n+1}\in \langle g_{1},\ldots,g_{n}\rangle \subset \langle y_{1},\ldots,y_{n}\rangle $. Since $\langle y_{1},\ldots,y_{n}\rangle$ is a prime ideal, we have $h\in \langle y_{1},\ldots,y_{n}\rangle$.\qed
\medskip

\begin{theorem}\label{gend}
Let $S = K[x_{ij}, \, y_{j} \mid 1\leq i\leq n+1,1\leq j\leq n]$ denote the polynomial 
$K$-algebra and $\widehat{X} = (x_{ij})_{(n+1)\times n}$, $Y=(y_{j})_{n\times 1}$ be 
generic matrices. Suppose $\mathcal{J} = I_{1}(\widehat{X}Y) = \langle g_{1},\ldots , g_{n+1}\rangle$, 
where $g_{i}=\sum_{j=1}^{n} x_{ij}y_{j}$. Let $\Delta_{i}$ denote the determinant 
of the matrix obtained by deleting $i$'th row of $\widehat{X}$ for $1\leq i\leq n+1$. 
\begin{enumerate}[(i)]
\item The sequence  $g_{1},\ldots ,g_{n+1}$ forms a $d$-sequence.
\item The Rees algebra of $I_{1}(\widehat{X}Y)$ is isomorphic to 
$S[T_{1},\ldots,T_{n+1}]/\mathfrak{I}$, where $\mathfrak{I}\subset S[T_{1},\ldots,T_{n}]$ 
is the ideal generated by the set 
$$\{T_{i}g_{j}-T_{j}g_{i},\displaystyle\sum_{i=1}^{n+1}-(1)^{k}\Delta_{k}T_{k}\mid 1\leq i,j\leq n+1\}.$$
\end{enumerate}
\end{theorem}

\proof (i) \, It is known by Theorem \ref{gensymn} that $g_{1},\ldots ,g_{n}$ is a regular sequence. 
Every regular sequence is also a $d$-sequence, therefore it is enough to prove that 
$$(\langle g_{0},\ldots,g_{i}\rangle:g_{i+1}g_{n+1})=(\langle g_{0},\ldots,g_{i}\rangle :g_{n+1}) $$ 
for all $0\leq i\leq n$, where $g_{0}=0$. By Lemma 6.7 in \cite{sstgrob}, we have, 
$(\langle g_{0},\ldots,g_{n}\rangle :g_{n+1})=\langle g_{1},\ldots,g_{n},\Delta_{n+1}\rangle$. 
Again by Theorem \ref{prime}, the ideals  $\langle g_{0},\ldots,g_{i}\rangle$ for $1\leq i\leq n-1$ 
and $\langle g_{1},\ldots,g_{n},\Delta_{n+1}\rangle$ are prime ideals. Therefore, we have 
$(\langle g_{0},\ldots,g_{i}\rangle:g_{i+1}g_{n+1})=(\langle g_{0},\ldots,g_{i}\rangle :g_{n+1})$ 
for all $0\leq i\leq n$. 
\medskip

\noindent (ii) \, A graded free minimal resolution of $I_{1}(\widehat{X}Y)$ can be found in 
Section 6 of \cite{sstgrob}, the construction 
of which uses the mapping cone technique. We now show, how this piece of information 
can be used to write the first syzygies matrix because no cancellation takes place 
at that level. We will be able to write the equations defining the Rees algebra if 
we can write the first syzygy matrix explicitly. 
\medskip

Let $P_{n}=\langle y_{1},\ldots,y_{n}\rangle$, $L_{n}=\langle g_{1},\ldots,g_{n}\rangle $, 
$T_{n}= \langle g_{1},\ldots,g_{n},\Delta_{n+1}\rangle$. Let $X_{ij}$ denote the 
$(ij)$-th cofactor of the matrix $X=(x_{ij})_{n\times n}$. We have the following 
diagram of exact chain complexes:
\medskip

{\scriptsize \xymatrix{
\cdots \ar[r] &(R(-n-2))^{\binom{n}{2}}\ar[r]^{\phi_{2n}}\ar[d]^{\delta_{2n}} & 
(R(-n-1))^{\binom{n}{1}}\ar[r]^{\phi_{1n}}\ar[d]^{\delta_{1n}} & R(-n)\ar[r]^{\phi_{0n}}\ar[d]^{\delta_{0n}= \Delta_{n}} & R/P_{n}\ar[d]\ar[r] & 0 \\
\cdots \ar[r] & (R(-4))^{\binom{n}{2}}\ar[r]_{d_{2}n} & (R(-2))^{\binom{n}{1}}\ar[r]_{d_{1}n} 
& R\ar[r]_{d_{0n}} & R/L_{n}\ar[r] & 0
}}

\noindent where $$\delta_{1n}=\begin{pmatrix}
X_{11} &\cdots & X_{1n}\\
\vdots &\vdots &\vdots\\
X_{n1} &\cdots & X_{nn}
\end{pmatrix},\,\, d_{1n}=\begin{pmatrix}
g_{1}&\cdots & g_{n}
\end{pmatrix},\,\, \phi_{1n}=\begin{pmatrix}
y_{1}&\cdots & y_{n}
\end{pmatrix}.
$$
By Lemma 6.5 in \cite{sstgrob}, the first box in the above diagram is commutative. 
It is evident that, for each $0\leq i\leq n$, the degrees of the domain and the 
co-domain of $\delta_{in}$ do not match. Therefore, a minimal graded free resolution 
of $T_{n}$ is obtained by the mapping cone technique:
$$
\cdots \longrightarrow (R(-4))^{\binom{n}{2}}\oplus (R(-n-1))^{\binom{n}{1}} \stackrel{\alpha_{2n}}\longrightarrow R(-2)^{n}\oplus R(-n) \stackrel{\alpha_{1n}}\longrightarrow R\longrightarrow R/T_{n}\longrightarrow 0$$
where $\alpha_{1n}= \begin{pmatrix}
g_{1} & \cdots & g_{n} & \Delta_{n}
\end{pmatrix}$ and $\, \alpha_{2n}= 
\left[
\begin{array}{c|c}
-\delta_{1n} & d_{2n} \\
\hline
\phi_{1n} & 0
\end{array}
\right]$. 
We now have the following diagram of exact chain complexes:
\medskip

{\scriptsize\xymatrix{
\cdots \ar[r]& (R(-6))^{\binom{n}{2}}\oplus (R(-n-3))^{\binom{n}{1}}\ar[r]^{\alpha_{2n}}\ar[d]^{\gamma_{2n}}&
R(-4)^{n}\oplus R(-n-2)\ar[r]^{\alpha_{1n}}\ar[d]^{\gamma_{1n}}&R(-2)\ar[r]^{\alpha_{0n}}\ar[d]^{\gamma_{0n}= g_{n+1}} &R/T_{n}\ar[d]\ar[r]&0 \\
\cdots \ar[r] &(R(-4))^{\binom{n}{2}}\ar[r]_{d_{2}n}&(R(-2))^{\binom{n}{1}}\ar[r]_{d_{1}n}&R\ar[r]_{d_{0n}}&R/L_{n}\ar[r]& 0
}}
\noindent where  \[
\gamma_{1n}=
\left[
\begin{array}{c|c}
-g_{n+1}I_{n} & B_{n} \\
\hline
d_{1n} & (-1)^{n+1}\Delta_{n+1}
\end{array}
\right], \quad \mbox{with} \quad B_{n}=\begin{pmatrix}
-\Delta_{1}\\
\vdots\\
(-1)^{n}\Delta_{n}
\end{pmatrix}.\] 
By Lemma \ref{cofac}, the first box is commutative in the above 
diagram. Therefore, a graded free resolution of $I_{1}(XY)$ is given by:
$$
\cdots \longrightarrow (R(-4))^{\binom{n}{2}+n}\oplus R(-n-2) \stackrel{\psi_{2n}}\longrightarrow R(-2)^{n+1} \stackrel{\psi_{1n}}\longrightarrow R\longrightarrow R/I_{1}(XY)\longrightarrow 0 
$$
where $\psi_{1n}= \begin{pmatrix}
g_{1}&\cdots &g_{n+1} 
\end{pmatrix}$ and $\psi_{2n}=
\left[
\begin{array}{c|c}
-\gamma_{1n} & d_{2n} \\
\hline
\alpha_{1n} & 0
\end{array}
\right]$.
\medskip

We know that the sequence $\{g_{1},\ldots,g_{n+1}\}$ is a $d$-sequence. Therefore, the 
ideal $I_{1}(\widehat{X}Y)$ is of linear type and the syzygy matrix $\psi_{2n}=
\left[
\begin{array}{c|c}
-\gamma_{1n} & d_{2n} \\
\hline
\alpha_{1n} & 0
\end{array}
\right]$ describes the ideal $\mathfrak{I}$ explicitly.\qed
\bigskip

\subsection{$I_{1}(XY)$ with alternating $X$.} Given an $n\times n$ generic alternating matrix $X$ and 
a generic column matrix $Y$, a scheme was proposed for computing a minimal free resolution of the 
ideal $I_{1}(XY)$ in the paper \cite{sstquad}. However, the scheme depended on two 
conjectures. We now prove those conjectures in Lemma \ref{altprime} and Lemma \ref{conj2} below. 
This not only validates our earlier work in \cite{sstquad} but also helps us prove Theorem \ref{altdseq} 
on the $d$-sequence. 
\medskip

We use the same 
notations as in \cite{sstquad}. Let $X_{n}$ denote the $n\times n$ alternating matrix and 
$Y_{n}$ denote the $n\times 1$ generic matrix given by 
$$X_{n}=\left[
\begin{array}{ccccc}
0& x_{12} &x_{13} & \ldots &x_{1n}\\
-x_{12}&0 & x_{23}& \ldots &x_{2n}\\
-x_{13}& -x_{23}&&\\
\vdots& \vdots &&\\
-x_{1n}& -x_{2n}&\ldots && 0
\end{array}
\right] \quad \mbox{and} \quad 
Y_{n}
 =\left[
\begin{array}{c}
y_{1}\\
\vdots\\
y_{n}
\end{array}
\right].$$
With the assumption that $x_{ij}=-x_{ji}$, if $i>j$ and $x_{ii}=0$, 
let $g_{ki}=\Sigma_{j=1}^{i}x_{kj}y_{j}$. Therefore 
$I_{1}(X_{n}Y_{n}) = \langle g_{1n},g_{2n},\cdots, g_{nn}\rangle$.  Let $\Delta_{(i)n}$ 
denote the Pffafian of the alternating matrix $X_{n}$ with 
the $i$-th row and the $i$-th column deleted.
\medskip

\begin{lemma}\label{useful} 
Let $x_{ij}=-x_{ji}$, if $i>j$ and $x_{ii}=0$. Then
\begin{itemize}
\item[(i)] $y_{n}g_{nn} = -\left( y_{1}g_{1n}+ y_{2}g_{2n}+\cdots+ y_{n-1}g_{(n-1)n}\right)$.

\item[(ii)] $g_{k (n-1)}g_{nn}  = x_{kn}y_{1}g_{1n}+ x_{kn}y_{2}g_{2n}+ \cdots + g_{nn}+x_{kn}y_{k}g_{kn}+\cdots + x_{kn}y_{n}g_{n-1 n}$.

\item[(iii)] $ \Delta_{(i)n}y_{i}  =  (-\Delta_{(1)n})g_{1n}+ (\Delta_{(2)n})g_{2n}+ \cdots +\widehat{(-1)^{i-1}\Delta_{(i)n})g_{in}}+\cdots ((-1)^{n-1}\Delta_{(n-1)n})g_{(n-1)n}$, \, for \, every \, $1\leq i\leq n$.
\end{itemize}
\end{lemma}

\proof These are easy to verify. \qed
\medskip

\begin{lemma}\label{akosul}
$\{g_{1n},\cdots , g_{(n-1) n} \}$ forms a regular sequence for $n\geq 2$.
\end{lemma}

\proof See part (ii) of Theorem 2.3 in \cite{sstprime}. \qed
\medskip

\begin{lemma}\label{altprime} Let $I_{n}=\langle g_{1n},g_{2n},\cdots , g_{(n-1) n}\rangle$ and 
$J_{n}= \langle g_{nn}\rangle$, so that $I_{n}+J_{n}=I_{1}(X_{n}Y_{n})$. 
Let 
$$C_{n}:= (I_{n}:J_{n})=\langle g_{1 (n-1)}, g_{2 (n-2)},\cdots, g_{n-1 (n-1)}, y_{n}, \Delta_{(n)n}\rangle .$$ 
If $n$ is even, then, $\Delta_{(n)n}=0$ and 
$C_{n}=\langle g_{1 (n-1)}, g_{2 (n-2)},\cdots, g_{(n-1) (n-1)}, y_{n}\rangle$, 
for every $n\geq 4$. Hence $C_{n}$ is a prime ideal for every $n\geq 4$.
\end{lemma}

\proof  At first we note that the ideal 
$$B_{n}:=\langle g_{1 (n-1)}, g_{2 (n-2)},\cdots, g_{n-1 (n-1)}, y_{n}, \Delta_{(n)n}\rangle$$ 
is a prime ideal by Proposition 5.8 and Lemma 5.12 in \cite{hu}. Let $hg_{nn}\in I_{n}\subset B_{n}$. 
We have $h\in B_{n}$, since $g_{nn}\notin B_{n}$. Thus $C_{n}\subset B_{n}$. The other inclusion 
easily follows from Lemma \ref{useful}. \qed
\medskip

\begin{lemma}\label{conj2} If $n$ is odd then $\Delta_{(n)n}\neq 0$, and for every $n\geq 4$,
$$P_{n}:= (\langle g_{1 (n-1)}, g_{2 (n-2)},\cdots, g_{(n-1) (n-1)}\rangle:\Delta_{(n)n})
= \langle y_{1},\cdots, y_{n-1} \rangle.$$
\end{lemma}

\proof Clearly $\langle y_{1},\cdots, y_{n-1} \rangle\subset P_{n}$ by the Lemma \ref{useful}. 
Let $h\Delta_{(n)n}\in I_{n-1}\subset \langle y_{1},\cdots, y_{n-1} \rangle$. Since 
$\Delta_{(n)n}\notin  \langle y_{1},\cdots, y_{n-1} \rangle$, we have $h\in \langle y_{1},\cdots, y_{n-1} \rangle $.\qed
\medskip

\begin{theorem}\label{altdseq}
$\{g_{1n},\cdots , g_{(n-1) n}, g_{nn}\}$ forms a  $d$-sequence for $n\geq 3$.
\end{theorem}

\proof We have $\{g_{1n},\cdots , g_{(n-1) n} \}$ forms a regular sequence for $n\geq 2$. 
Therefore, it is enough to show that, for $1\leq i\leq n-1$ and $j\geq i+1$,
$$(\langle g_{0n},g_{1n},\cdots , g_{in}\rangle:g_{(j+1)n}g_{nn} )=
\langle(g_{0n},g_{1n},\cdots , g_{in}\rangle:g_{nn} ),$$ 
where $g_{0n}=0$ for all $n$. We have, for every $j\geq i+1$,
$$(\langle g_{0n},g_{1n},\cdots , g_{in}\rangle:g_{nn} )\subseteq (\langle g_{0n},g_{1n},\cdots , g_{in}\rangle:g_{(j+1)n}g_{nn} ).$$ 
Let $h\in (\langle g_{0n},g_{1n},\cdots , g_{in}\rangle:g_{(j+1)n}g_{nn} )$, then $hg_{(j+1)n}g_{nn}\in \langle g_{0n},g_{1n},\cdots , g_{in}\rangle$. For $1\leq i\leq n-2$, the ideal $\langle(g_{0n},g_{1n},\cdots , g_{in}\rangle$ 
is a prime ideal (see Theorem 3.3 in \cite{sstprime}), therefore $hg_{nn}\in \langle g_{0n},g_{1n},\cdots , g_{in}\rangle$. 
Moreover, for $i=n-1$, the ideal $(\langle g_{0n},g_{1n},\cdots , g_{(n-1)n}\rangle:g_{nn} )$ is a prime ideal by 
Lemma \ref{altprime}, therefore $h\in (\langle g_{0n},g_{1n},\cdots , g_{(n-1)n}\rangle:g_{nn} )$. \qed
\medskip

As an upshot of the aforesaid results, we get a primary decomposition of the ideal $I_{1}(X_{n}Y_{n})$ and 
somemore information regarding its normality.
\medskip 

\begin{theorem}\label{primary} The primary decomposition of $$I_{1}(X_{n}Y_{n})=\langle g_{1n},\ldots,g_{nn},\Delta_{n}\rangle\cap \langle y_{1},\ldots,y_{n}\rangle$$ where $\Delta_{n}$ is the pfaffian  of the matrix $X_{n}$.
\end{theorem}
\proof Let $\mathfrak{Q}_{1n}=\langle g_{1n},\ldots,g_{nn},\Delta_{n}\rangle$ and $\mathfrak{Q}_{2n}=\langle y_{1},\ldots,y_{n} \rangle$. We note that $\mathfrak{Q}_{2n}$ is a prime ideal and by proposition 5.8 and lemma 5.12 in \cite{hu}, $\mathfrak{Q}_{1n}$ is a prime ideal. Obviously $I_{1}(X_{n}Y_{n}) \subset  \mathfrak{Q}_{1n}\cap \mathfrak{Q}_{2n}$. Let $h\in \mathfrak{Q}_{1n}\cap \mathfrak{Q}_{2n}$,  where $h=\displaystyle \sum_{k=1}^{n}a_{k}g_{kn}+b\Delta_{n}$. Then $b\Delta_{n}\in \mathfrak{Q}_{2n}$. Hence $b\in \mathfrak{Q}_{2n}$. By lemma \ref{useful} $b\Delta_{n}\in I_{1}(X_{n}Y_{n})$, therefore $h\in I_{1}(X_{n}Y_{n})$. \qed
\medskip

\begin{corollary}
The ideal $I_{1}(X_{n}Y_{n})$ is a radical ideal.
\end{corollary}

\proof Follows from Theorem \ref{primary}.\qed
\medskip

\begin{definition}
Let $R$ be a Noetherian ring and $I$ an ideal. Then $I$ is normally torsionfree if 
$\mbox{Ass} R/I = \mbox{Ass} R/I^{n}$ for $n\geq 1$.
\end{definition}

\begin{theorem}\label{normal}
Let $R$ be a regular local ring and $I$ a reduced ideal. If $I$ is normally torsionfree, then $I$ is
normal.
\end{theorem}

\proof See the Proposition 1.54 in \cite{vas}.\qed
\medskip

\begin{lemma}
The ideal $I_{1}(X_{n}Y_{n})$ is normal.
\end{lemma}

\proof Follows from Theorems \ref{primary} and \ref{normal}.\qed
\bigskip

\section{Construction of a regular sequence using $I_{1}(XY)$}

\begin{lemma}\label{disjoint}
Let $h_{1},h_{2}\cdots, h_{n}\in R$ be such that with respect to a suitable 
monomial order on $R$ the leading terms of them are mutually coprime. Then, 
$h_{1},h_{2}\cdots, h_{n}$ is a regular sequence in $R$.
\end{lemma}

\proof. See lemma 2.2 in \cite{sstprime}.\qed
\medskip

\begin{lemma}\label{technical}
Let $h_{1},\ldots,h_{n-1}$ be distinct polynomials and for 
$1\leq r\leq n-1$ let $m_{1},\ldots m_{r+1}$ be distinct monomials 
in $R$. Suppose that the following properties are satisfied with 
respect to some monomial order on $R$:

\begin{enumerate}
\item[(i)] ${\rm Lt}(h_{i})=a_{i}b_{i}$ for every $1\leq i\leq n-1$; 

\item[(ii)] $\gcd({\rm Lt}(h_{i}), {\rm Lt}(h_{j})) = 1$ for every $1\leq i\neq j\leq n-1$;

\item[(iii)] $\gcd(m_{i}, m_{j}) = 1$ \, for every $1\leq i\neq j\leq r+1$;

\item[(iv)] $\gcd({\rm Lt}(h_{i}), m_{j}) = 1$ \, for every $1\leq i\leq n-1$ and 
$1\leq j\leq r+1$.
\end{enumerate}
Let $h_{n}=b_{k_{1}}m_{1}+\ldots b_{k_{r}}m_{r}+m_{r+1}+(\textrm{lower order terms})$. 
Then 

\begin{enumerate}
\item $h_{1},\ldots,h_{n-1},m_{1},\ldots m_{r},h_{n}$ is a regular sequence.  

\item Moreover, if $g$ is a polynomial such that $\gcd(\LT(g),h_{i})=1$ for 
$1\leq i\leq n-1$ and $\gcd(\LT(g),m_{i})=1$ for 
$1\leq i\leq r+1$, then $h_{1},\ldots,h_{n-1},m_{1},\ldots , m_{r},h_{n},g$ is 
a regular sequence.
\end{enumerate}
\end{lemma}

\proof \textbf{(1)} The sequence $h_{1},\ldots,h_{n-1},m_{1},\ldots ,m_{r}$ is 
a regular sequence by the Lemma \ref{disjoint}. Let 
$\widetilde{h}_{n}=h_{n}-\sum_{i=1}^{r}b_{k_{i}}m_{i}$, then 
$\LT(\tilde{h}_{n})=m_{r+1}$ is coprime with $\LT(h_{1}),\ldots,\LT(h_{n-1})$ and also coprime with $m_{1},\ldots , m_{r}$. Again by the Lemma \ref{disjoint}, we can write 
$h_{1},\ldots,h_{n-1},m_{1},\ldots ,m_{r},\widetilde{h}_{n}$ is a regular sequence.
\medskip

Let $h_{n}\cdot p=\sum_{i=1}^{n-1}h_{i}p_{i}+\sum_{i=1}^{r}m_{i}q_{i}$, therefore 
$$\widetilde{h}_{n}\cdot p=\sum_{i=1}^{n-1}h_{i}p_{i}
+\sum_{i=1}^{r}m_{i}(q_{i}-b_{k_{i}}).$$ 
This gives $\, p\in \langle h_{1},\ldots,h_{n-1},m_{1},\ldots ,m_{r}\rangle$, 
since $h_{1},\ldots,h_{n-1},m_{1},\ldots ,m_{r},\widetilde{h}_{n}$ is a regular 
sequence. Therefore, $h_{1},\ldots,h_{n-1},m_{1},\ldots m_{r},h_{n}$ is a regular 
sequence. 
\medskip

\noindent\textbf{(2)} Let $g\cdot p'=\sum_{i=1}^{n}h_{i}p'_{i}+\sum_{i=1}^{r}m_{i}q'_{i}$, therefore 
$$g\cdot p'=\sum_{i=1}^{n-1}h_{i}p'_{i}+\sum_{i=1}^{r}m_{k_{i}}(q'_{i}+b_{k_{i}})+\widetilde{h}_{n}p'_{n}.$$ 
Now  $h_{1},\ldots,h_{n-1},m_{1},\ldots m_{r},\widetilde{h}_{n},g$ \, is a regular 
sequence by \ref{disjoint}; hence $p'\in \langle h_{1},\ldots,h_{n-1},m_{k_{1}},\ldots ,m_{r},\tilde{h}_{n}\rangle$.\qed
\medskip

\begin{theorem}\label{maintheorem}
Suppose that $k_{t}=1+\left(\left[\frac{n}{t}\right]-1\right)t$, then $k_{n}=1$. 
The ordered set $$\mathcal{F} = \lbrace f_{11},\ldots , f_{n1}\rbrace \cup 
\lbrace f_{12},f_{32},\ldots, f_{k_{2}2}\rbrace \cup  
\lbrace f_{13}, f_{43}, \ldots , f_{k_{3}3}\rbrace \cup \cdots \cup 
\lbrace f_{1n}\rbrace$$ 
is a regular sequence in $R$. The polynomials occurring in the list follow 
the pattern indicated below:
$$\begin{pmatrix}
f_{11} & f_{12} & f_{13} & f_{14} & \cdots & f_{1n}\\
f_{21} & \times & \times & \times & \cdots & \times\\
f_{31} & f_{32} & \times & \times & \cdots & \times\\
f_{41} & \times & f_{43} & \times & \cdots & \times\\
f_{51} & f_{52} & \times & f_{54} & \cdots & \times\\
f_{61} & \times & \times & \times & \cdots & \times\\
f_{71} & f_{72} & f_{73} & \times & \cdots & \times\\
f_{81} & \times & \times & \times & \cdots & \times\\
\vdots & \vdots & \vdots & \vdots & \cdots & \vdots
\end{pmatrix}$$
\end{theorem}

\proof Consider the lexicographic monomial order given by 
\begin{eqnarray*}
x_{11} \cdots >x_{nn} & > & x_{12}>x_{23}>\cdots > x_{(n-1)n}\\
{} &  > & x_{13}>x_{24} > \cdots > x_{(n-2)n}\\
{} & > & \vdots\\
{} & > & x_{1n}\\
{} & > & x_{ij} > y_{kl}, \quad \text{for \,all} \quad i>j\quad \text{and}\quad 1\leq k,l \leq n.
\end{eqnarray*}

\noindent In order to show that the set $\mathcal{F}$ is a regular sequence we consider a 
larger ordered set of polynomials by adding some indeterminates in the list 
so that they follow the properties listed in Lemma \ref{technical}. The new 
ordered set of polynomials we consider is  
\begin{eqnarray*}
\widetilde{\mathcal{F}} & = & \{f_{11},\ldots f_{n1}\}\cup\\
{} & {} & \{y_{12},f_{12},y_{32},f_{32},\ldots,y_{k_{2}2},f_{k_{2}2}\}\cup\\
{} & {} & \{y_{13},y_{23},f_{13},\ldots,y_{k_{3}3},y_{(k_{3}+1)3},f_{k_{3}3}\}\cup\\ 
{} & {} & \vdots\\
{} & {} & \{y_{1t},\ldots,y_{(t-1)t}, f_{1t},\ldots ,f_{k_{t}t}\}\cup \\
{} & {} & \vdots\\
{} & {} & \{y_{1n},\ldots,y_{(n-1)n}, f_{1n}\}.
\end{eqnarray*}
If we can show that $\widetilde{\mathcal{F}}$ is a regular sequence in this order, 
then, under any permutation of this order the polynomials would still form a regular sequence because of homogeneity. Therefore, we can rearrange the entries of 
$\widetilde{\mathcal{F}}$ in such a way that the elements that appear first are the 
ones listed in $\mathcal{F}$ and that proves our claim. By Lemma \ref{disjoint}, 
$f_{11},\ldots f_{n1},y_{12}$ is a regular sequence. Using Lemma \ref{technical} 
we can add $f_{12}=x_{11}y_{12}+x_{12}y_{22}+\sum_{k=3}^{n}x_{1k}y_{k2}$ in the 
list. Therefore, $f_{11},\ldots f_{n1},y_{12}, f_{12}$ is a regular sequence. 
The indeterminate $y_{32}$ does not divide the any of the leading terms of 
$f_{11},\ldots f_{n1},y_{12},f_{12}-x_{11}y_{12}$. Therefore, $f_{11},\ldots f_{n1},y_{12},f_{12},y_{32}$ is a regular sequence, by Lemma  
\ref{technical}. Therefore, $\widetilde{\mathcal{F}}$ is a regular 
sequence if we continue the same process. The above comment proves that $\mathcal{F}$ is a regular sequence.\qed

\bibliographystyle{amsalpha}

\begin{thebibliography}{A}

\bibitem{herzog}
{J. Herzog, \emph{Certain Complexes Associated to a Sequence and a Matrix}, 
Manuscipta Math. 12(1974) 217--248.
}

\bibitem{hu}
{C. Huneke and B. Ulrich, \emph{Divisor Class Groups and Deformations}, 
American Journal of Mathematics 107(6)(1985)1265--1303.
} 

\bibitem{sstsum}
{J. Saha, I. Sengupta, G. Tripathi, \emph{Transversal Intersection and Sum of Polynomial Ideals}, 
arXiv:1611.04732v2 [math.AC] 2018.
}

\bibitem{sstprime} 
{J. Saha, I. Sengupta, G. Tripathi, \emph{Primary decomposition and normality
of certain determinantal ideals}, Proc. Indian Acad. Sci. (Math. Sci.) 129:55(2019). 
}

\bibitem{sstgrob}
{J. Saha, I. Sengupta, G. Tripathi, \emph{Ideals of the form $I_{1}(XY)$}, 
Journal of Symbolic Computation 91(2019)17--29.
}

\bibitem{sstquad}
{J. Saha, I. Sengupta, G. Tripathi, \emph{Quadrics defined by Skew-Symmetric Matrices}, 
International Journal of Algebra 11(8)(2017) 349 -- 356.
}

\bibitem{hunekeswanson}
{I. Swanson, C. Huneke, \emph{Integral Closure of Ideals, Rings, and Modules}, LMS 
Lecture Note Series 336. Cambridge University Press, UK, 2006.
}

\bibitem{vas}
{W. Vasconcelos, \emph{Integral Closure}, Springer Monograph in Mathematics, 2005.
}
\end{thebibliography}

\end{document}